\documentclass[11pt]{article}
\usepackage{amsmath,amssymb}
\newtheorem{theo}{Theorem}[section]

\newtheorem{lemma}[theo]{Lemma}

\newcommand{\ignore}[1]{}
\def\square{\vrule height6pt width7pt depth1pt}
\def\endpf{\hfill\square\bigskip}

\begin{document}
\title{Integer and fractional packing of families of graphs}
\author{Raphael Yuster
\thanks{
e-mail: raphy@research.haifa.ac.il \qquad
World Wide Web: http:$\backslash\backslash$research.haifa.ac.il$\backslash$\symbol{126}raphy}
\\ Department of Mathematics\\ University of
Haifa at Oranim\\ Tivon 36006, Israel}

\date{} 

\maketitle
\setcounter{page}{1}
\begin{abstract}
Let ${\cal F}$ be a family of graphs. For a graph $G$, the {\em ${\cal F}$-packing number}, denoted
$\nu_{{\cal F}}(G)$, is the maximum number of pairwise
edge-disjoint elements of ${\cal F}$ in $G$. A function $\psi$ from the set of elements of ${\cal F}$ in $G$ to $[0,1]$
is a {\em fractional ${\cal F}$-packing} of $G$ if $\sum_{e \in H \in {\cal F}} {\psi(H)} \leq 1$ for each $e \in E(G)$.
The {\em fractional ${\cal F}$-packing number}, denoted $\nu^*_{{\cal F}}(G)$, is defined to be the maximum value
of $\sum_{H \in {{G} \choose {{\cal F}}}} \psi(H)$ over all fractional ${\cal F}$-packings $\psi$.
Our main result is that $\nu^*_{{\cal F}}(G)-\nu_{{\cal F}}(G) = o(|V(G)|^2)$.
Furthermore, a set of $\nu_{{\cal F}}(G) -o(|V(G)|^2)$ edge-disjoint elements of ${\cal F}$ in $G$ can be found
in randomized polynomial time.
For the special case ${\cal F}=\{H_0\}$ we obtain a significantly simpler proof of a recent difficult result of Haxell and R\"odl \cite{HaRo}
that $\nu^*_{H_0}(G)-\nu_{H_0}(G) = o(|V(G)|^2)$.
\end{abstract}

\section{Introduction}
All graphs considered here are finite and have no loops, multiple edges or isolated vertices.
For the standard terminology used the reader is referred to \cite{Bo}.
Let ${\cal F}$ be any fixed finite or infinite family of graphs. For a graph $G$, the {\em ${\cal F}$-packing number}, denoted
$\nu_{{\cal F}}(G)$, is the maximum number of pairwise
edge-disjoint copies of elements of ${\cal F}$ in $G$. A function $\psi$ from the set of copies of elements of ${\cal F}$ in $G$ to $[0,1]$
is a {\em fractional ${\cal F}$-packing} of $G$ if $\sum_{e \in H \in {\cal F}} {\psi(H)} \leq 1$ for each $e \in E(G)$.
For a fractional ${\cal F}$-packing $\psi$, let $w(\psi)=\sum_{H \in {{G} \choose {{\cal F}}}} \psi(H)$.
The {\em fractional ${\cal F}$-packing number}, denoted $\nu^*_{{\cal F}}(G)$, is defined to be the maximum value
of $w(\psi)$ over all fractional packings $\psi$. Notice that, trivially, $\nu^*_{{\cal F}}(G) \geq \nu_{{\cal F}}(G)$.
If ${\cal F}$ consists of a single graph $H_0$ we shall denote the parameters above by $\nu_{H_0}(G)$ and $\nu^*_{H_0}(G)$.

Since computing $\nu^*_{{\cal F}}(G)$ amounts to solving a linear program, it can be computed in polynomial time for every finite
${\cal F}$. On the other hand, it was proved by Dor and Tarsi in \cite{DoTa} that computing $\nu_{H_0}(G)$ is NP-Hard for every
$H_0$ with a component having at least three edges.
Thus, it is interesting to determine when $\nu^*_{{\cal F}}(G)$ and $\nu_{{\cal F}}(G)$ are ``close'',
thereby getting a polynomial time approximating algorithm for an NP-Hard problem. The following result was proved by Haxell and R\"odl in \cite{HaRo}.
\begin{theo}
\label{t0}
If $H_0$ is a fixed graph and $G$ is a graph with $n$ vertices, then $\nu^*_{H_0}(G)-\nu_{H_0}(G) = o(n^2)$.
\end{theo}
The 25 page proof of Theorem \ref{t0} presented in \cite{HaRo} is very difficult.
The major difficulty lies in the fact that their method requires proving that there is a fractional packing which is only slightly less than optimal,
and which assigns to every copy of $H_0$ either $0$ or a value greater than $\tau$ for some $\tau > 0$ which is only a function of $H_0$.

In this paper we present a significantly simpler proof of Theorem \ref{t0}. Our proof method enables us to generalize Theorem \ref{t0} to
the ``family'' case. There does not seem to be an easy way to generalize the proof in \cite{HaRo} to the family case.
\begin{theo}
\label{t1}
If ${\cal F}$ is a fixed family of graphs and $G$ is a graph with $n$ vertices, then $\nu^*_{{\cal F}}(G)-\nu_{{\cal F}}(G) = o(n^2)$.
\end{theo}
Notice that Theorem \ref{t1} immediately yields a polynomial time algorithm for approximating $\nu_{{\cal F}}(G)$ to within
an additive term of $\epsilon n^2$ for every $\epsilon > 0$. Furthermore, if ${\cal F}$ is finite, the degree of the polynomial depends only
on ${\cal F}$, and not on $1/\epsilon$. Our proof also supplies a randomized polynomial time algorithm that {\em finds} a set of
$\nu_{{\cal F}}(G)-o(n^2)$ edge-disjoint copies of elements of ${\cal F}$ in $G$.

\section{Tools used in the main result}
As in \cite{HaRo}, a central ingredient in our proof of the main result is Szemer\'edi's regularity
lemma \cite{Sz}. Let $G=(V,E)$ be a graph, and let $A$ and $B$ be two disjoint subsets
of $V(G)$. If $A$ and $B$ are non-empty, let $E(A,B)$ denote set of edges between them, and
put $e(A,B)=|E(A,B)|$. The {\em density of edges} between $A$ and $B$ is defined as
$$
d(A,B) = \frac{e(A,B)}{|A||B|}.
$$
For $\gamma>0$ the pair $(A,B)$ is called {\em $\gamma$-regular}
if for every $X \subset A$ and $Y \subset B$ satisfying
$|X|>\gamma |A|$ and $|Y|>\gamma |B|$ we have
$$
|d(X,Y)-d(A,B)| < \gamma.
$$
An {\em equitable partition} of a set $V$ is a partition of $V$ into
pairwise disjoint classes $V_1,\ldots,V_m$ whose sizes are as equal as possible.
An equitable partition
of the set of vertices $V$ of a graph $G$ into the classes $V_1,\ldots,V_m$ is
called {\em $\gamma$-regular} if $|V_i| < \gamma |V|$ for
every $i$ and all but at most $\gamma {m \choose 2}$ of the pairs
$(V_i,V_j)$ are $\gamma$-regular.
The regularity lemma states the following:
\begin{lemma}
\label{l21}
For every $\gamma>0$, there is an integer $M(\gamma)>0$ such that for
every graph $G$ of order $n > M$ there is a $\gamma$-regular partition of
the vertex set of $G$ into $m$ classes, for some $1/\gamma < m < M$. \endpf
\end{lemma}

Let $H_0$ be a fixed graph with the vertices $\{1,\ldots,k\}$, $k \geq 3$.
Let $W$ be a $k$-partite graph with vertex classes $V_1,\ldots,V_k$. A subgraph
$J$ of $W$ with ordered vertex set $v_1,\ldots,v_k$ is {\em partite-isomorphic} to $H_0$ if $v_i \in V_i$
and the map $v_i \rightarrow i$ is an isomorphism from $J$ to $H_0$.

The following lemma is almost identical to the (2 page) proof of Lemma 15 in \cite{HaRo}
and hence the proof is omitted.
\begin{lemma}
\label{l22}
Let $\delta$ and $\zeta$ be positive reals.
There exist $\gamma=\gamma(\delta, \zeta, k)$ and
$T=T(\delta, \zeta, k)$ such that the following holds.
Let $W$ be a $k$-partite graph with vertex classes $V_1,\ldots,V_k$
and $|V_i|=t > T$ for $i=1,\ldots,k$. Furthermore, for each $(i,j) \in E(H_0)$, $(V_i,V_j)$ is
a $\gamma$-regular pair with density $d(i,j) \geq \delta$ and for each $(i,j) \notin E(H_0)$,
$E(V_i,V_j)=\emptyset$.
Then, there exists a spanning subgraph $W'$ of $W$, consisting of at least $(1-\zeta)|E(W)|$ edges
such that the following holds.
For an edge $e \in E(W')$, let $c(e)$ denote the number of subgraphs of $W'$ that are partite isomorphic
to $H_0$ and that contain $e$.
Then, for all $e \in E(W')$, if $e \in E(V_i,V_j)$ then
$$
\left|c(e) - t^{k-2} \frac{\prod_{(s,p) \in E(H_0)}d(s,p)}{d(i,j)}\right| < \zeta t^{k-2}.
$$
\endpf
\end{lemma}

Finally, we need to state the seminal result of Frankl and R\"odl \cite{FrRo} on
near perfect coverings and matchings of uniform hypergraphs. Recall that if $x,y$ are two vertices
of a hypergraph then $deg(x)$ denotes the degree of $x$ and $deg(x,y)$ denotes the number of edges  that contain
both $x$ and $y$ (their {\em co-degree}). We use the version of the Frankl  and R\"odl Theorem due to
Pippenger (see, e.g., \cite{Fu}).
\begin{lemma}
\label{l23}
For an integer $r \geq 2$ and a real $\beta > 0$ there exists a real $\mu > 0$ so that:
If the $r$-uniform hypergraph $L$ on $q$ vertices has the following properties for some $d$:\\
(i) $(1-\mu)d < deg(x) < (1+\mu)d$ holds for all vertices,\\
(ii) $deg(x,y) < \mu d$ for all distinct $x$ and $y$,\\
then $L$ has a matching of size at least $(q/r)(1-\beta)$. \endpf
\end{lemma}

\section{Proof of the main result}
Let ${\cal F}$ be a family of graphs, and let $\epsilon > 0$.
To avoid the trivial case we assume $K_2 \notin {\cal F}$.
We shall prove there exists $N=N({\cal F},\epsilon)$ such that for all
$n > N$, if $G$ is an $n$-vertex graph then $\nu_{{\cal F}}^*(G) - \nu_{{\cal F}}(G) < \epsilon n^2$.

Let $k_\infty$ denote the maximal order of a graph in ${\cal F}$. Let $k_0=\min\{k_\infty,\lceil 20/\epsilon \rceil\}$.
Let $\delta=\beta=\epsilon/4$. For all $r=2,\ldots,{k_0}^2$, let $\mu_r=\mu(\beta,r)$ be as in Lemma \ref{l23}, and put
$\mu=\min_{r=2}^{k_0^2} \{\mu_r\}$.
Let $\zeta=\mu \delta^{{k_0}^2}/2$. For $k=3, \ldots, k_0$, let $\gamma_k=\gamma(\delta, \zeta, k)$ and
$T_k=T(\delta, \zeta, k)$ be as in Lemma \ref{l22}. Let $\gamma=\min_{k=3}^{k_0} \{\gamma_k\}$.
Let $M=M(\gamma\epsilon/(25{k_0}^2))$ be as in Lemma \ref{l21}.
Finally, we shall define $N$ to be a sufficiently large constant, depending on the above chosen parameters,
and for which various conditions stated in the proof below hold (it will be obvious in the proof that all these conditions
indeed hold for $N$ sufficiently large). Thus, indeed, $N=N({\cal F}, \epsilon)$.

Fix an $n$-vertex graph $G$ with $n > N$ vertices. Fix a fractional ${\cal F}$-packing
$\psi$ with $w(\psi)=\nu_{{\cal F}}^*(G)$. We may assume that $\psi$ assigns a
value to each {\em labeled} copy of an element of ${\cal F}$ simply by dividing the value of $\psi$
on each nonlabeled copy by the size of the automorphism group of that element.
If $\nu_{{\cal F}}^*(G) < \epsilon n^2$ we are done.
Hence, we assume $\nu_{{\cal F}}^*(G) =\alpha n^2 \geq \epsilon n^2$.

We apply Lemma \ref{l21} to $G$ and obtain a $\gamma'$-regular
partition with $m'$ parts, where $\gamma'=\gamma \epsilon/(25{k_0}^2)$ and $1/\gamma' < m' < M(\gamma')$.
Denote the parts by $U_1,\ldots,U_{m'}$.
Notice that the size of each part is either $\lfloor n/{m'} \rfloor$ or $\lceil n/{m'} \rceil$.
For simplicity we may and will assume that $n/{m'}$ is an integer, as this assumption does not
affect the asymptotic nature of our result. For the same reason we may and will assume that
$n/(25m'{k_0}^2/\epsilon)$ is an integer.

We randomly partition each $U_i$ into $25{k_0}^2/\epsilon$ equal parts of size $n/(25m'{k_0}^2/\epsilon)$ each.
All $m'$ partitions are independent. We now have $m=25m'{k_0}^2/\epsilon$ {\em refined} vertex classes, denoted
$V_1,\ldots,V_m$.
Suppose $V_i \subset U_s$ and $V_j \subset U_t$ where $s \neq t$. We claim that if
$(U_s,U_t)$ is a $\gamma'$-regular pair then $(V_i,V_j)$ is a $\gamma$-regular pair.
Indeed, if $X \subset V_i$ and $Y \subset V_j$ have $|X|, |Y| > \gamma n/(25m'{k_0}^2/\epsilon)$ then
$|X|, |Y| > \gamma' n/m'$ and so $|d(X,Y) - d(U_s,U_t)| < \gamma'$. Also
$|d(V_i,V_j) - d(U_s,U_t)| < \gamma'$. Thus, $|d(X,Y) - d(V_i,V_j)| < 2\gamma' < \gamma$.

Let $H$ be a labeled copy of some $H_0 \in {\cal F}$ in $G$. If $H$ has $k$ vertices and $k \leq k_0$ then
the expected number of pairs of vertices of $H$
that belong to the same vertex class in the refined partition is clearly at most
${k \choose 2}\epsilon/(25{k_0}^2) < \epsilon/50$. Thus, the probability that $H$ has two vertices in the
same vertex class is also at most $\epsilon/50$. We call $H$ {\em good} if it has $k \leq k_0$ vertices and its $k$ vertices belong to
$k$ distinct vertex classes of the refined partition.
By the definition of $k_0$, if $H$ has $k > k_0$ vertices and $\psi(H) > 0$ then we must have $k > 20/\epsilon$.
Since graphs with $k$ vertices have at least $k/2$ edges, the contribution of graphs with $k > k_0$ vertices to
$\nu_{{\cal F}}^*(G)$ is at most ${n \choose 2}/(10/\epsilon) < \epsilon n^2/20$.
Hence, if $\psi^{**}$ is the restriction of $\psi$ to
good copies (the bad copies having $\psi^{**}(H)=0$) then the expectation of $w(\psi^{**})$ is at least
$(\alpha-\epsilon/50-\epsilon/20)n^2$. We therefore {\em fix} a partition $V_1,\ldots,V_m$ for which
$w(\psi^{**}) \geq (\alpha-0.07\epsilon)n^2$.

Let $G^*$ be the spanning subgraph of $G$ consisting of the edges with endpoints
in distinct vertex classes of the refined partition that form a $\gamma$-regular pair with density at least $\delta$
(thus, we discard edges inside classes, between non regular pairs, or between sparse pairs).
Let $\psi^*$ be the restriction of $\psi^{**}$ to the labeled copies of elements of ${\cal F}$ in $G^*$.
We claim that
$\nu_{{\cal F}}^*(G^*) \geq w(\psi^*) > w(\psi^{**})-0.72\delta n^2 \geq (\alpha-0.07\epsilon-0.72\delta)n^2=(\alpha - \delta)n^2$.
Indeed, by considering the number of discarded edges we get (using $m' > 1/\gamma'$ and $\delta >> \gamma'$)
$$
w(\psi^{**}) - w(\psi^*) \leq |E(G) - E(G^*)|  <
\gamma' {{m'} \choose 2}\frac{n^2}{{m'}^2} +
{{m'} \choose 2}(\delta+\gamma')\frac{n^2}{{m'}^2} +
{m'}{{n/{m'}} \choose 2} < 0.72 \delta n^2.
$$
Let $R$ denote the $m$-vertex graph whose vertices are $\{1,\ldots,m\}$ and $(i,j) \in E(R)$
if and only if $(V_i,V_j)$ is a $\gamma$-regular pair with density at least $\delta$.
We define a (labeled) fractional ${\cal F}$-packing $\psi'$ of $R$ as follows.
Let $H$ be a labeled copy of some $H_0 \in {\cal F}$ in $R$ and assume that the vertices of $H$ are
$\{u_1,\ldots,u_k\}$ where $u_i$ plays the role of vertex $i$ in $H_0$.
We define $\psi'(H)$ to be the sum of the values of $\psi^*$ taken over all
subgraphs of $G^*[V_{u_1},\ldots,V_{u_k}]$ which are partite isomorphic to $H_0$,
divided by $n^2/m^2$. Notice that by normalizing with $n^2/m^2$ we guarantee
that $\psi'$ is a proper fractional ${\cal F}$-packing of $R$ and that
$\nu_{{\cal F}}^*(R) \geq w(\psi') =m^2w(\psi^*)/n^2 \geq m^2(\alpha - \delta)$.

We use $\psi'$ to define a random coloring of the edges of $G^*$.
Our ``colors'' are the labeled copies of elements of ${\cal F}$ in $R$.
Let $d(i,j)$ denote the density of $(V_i,V_j)$ and
notice that $|E_{G^*}(V_i,V_j)|=d(i,j)n^2/m^2$. Let $H$ be a labeled copy of some $H_0 \in {\cal F}$ in
$R$, and assume that $H$ contains the edge $(i,j)$. Each $e \in E(V_i,V_j)$ is chosen to
have the ``color'' $H$ with probability $\psi'(H) /d(i,j)$.
The choices made by distinct edges of $G^*$ are independent.
Notice that this random coloring is legal (in the sense that the sum of probabilities is
at most one) since the sum of $\psi'(H)$ taken
over all labeled copies of elements of ${\cal F}$ containing $(i,j)$ is at most $d(i,j) \leq 1$.
Notice also that some edges might stay uncolored in our random coloring of the edges of $G^*$.

Let $H$ be a labeled copy of some $H_0 \in {\cal F}$ in $R$, and assume that $\psi'(H) > m^{1-k_0}$.
Without loss of generality, assume that the vertices of $H$ are $\{1,\ldots,k\}$ where $i \in V(H)$ plays
the role of $i \in V(H_0)$. Let $r$ denote the number of edges of $H$. Notice that $r < k_0^2$.
Let $W_H=G^*[V_1,\ldots,V_k]$ (in fact we only consider edges between pairs that correspond to edges of $H_0$).
Notice that $W_H$ is a subgraph of $G^*$ which satisfies the conditions in Lemma \ref{l22},
since $t=n/m > N\epsilon/(25{k_0}^2M) > T_k$ (here we assume $N > 25{k_0}^2MT_k/\epsilon$).
Let $W'_H$ be the spanning subgraph of $W_H$ whose existence is guaranteed in Lemma \ref{l22}.
Let $X_H$ denote the spanning subgraph of $W'_H$ consisting only of the edges whose color is $H$.
Notice that $X_H$ is a random subgraph of $W'_H$.
For an edge $e \in E(X_H)$, let $C_H(e)$ denote the set of subgraphs of $X_H$ that contain $e$
and that are partite isomorphic to $H_0$. Put $c_H(e)=|C_H(e)|$. A crucial argument is the following:
\begin{lemma}
With probability at least $1-m^3/n$, for all $e \in E(X_H)$,
\label{l31}
\begin{equation}
\label{e0}
\left|c_H(e) - t^{k-2} \psi'(H)^{r-1} \right| < \mu \psi'(H)^{r-1} t^{k-2}.
\end{equation}
\end{lemma}
{\bf Proof:}\,
Let $C(e)$ denote the set of subgraphs of $W'_H$ that contain $e$ and that are partite isomorphic to $H_0$.
Put $c(e)=|C(e)|$. According to Lemma \ref{l22}, if $e \in E(V_i,V_j)$ then
\begin{equation}
\label{e1}
\left|c(e) - t^{k-2} \frac{\prod_{(s,p) \in E(H_0)}d(s,p)}{d(i,j)}\right| < \zeta t^{k-2}.
\end{equation}
Fix an edge $e \in E(X_H)$ belonging to $E(V_i,V_j)$.
The probability that an element of $C(e)$ also belongs to $C_H(e)$ is precisely
$$
\rho=\psi'(H)^{r-1} \cdot \frac{d(i,j)}{\prod_{(s,p) \in E(H_0)}d(s,p)}.
$$
We say that two distinct elements $Y,Z \in C(e)$ are {\em dependent} if they share at least one edge
other than $e$.
Consider the dependency graph $B$ whose vertex set is $C(e)$ and the edges connect dependent pairs.
Since two dependent elements share at least three vertices (including the two endpoints of $e$), we have
that $\Delta(B) = O(t^{k-3})$. Hence, $\chi(B)=O(t^{k-3})$. Put $s=\chi(B)$. Let $C^1(e), \ldots, C^s(e)$ denote
a partition of $C(e)$ to independent sets. Let $C^q_H(e) =C^q(e) \cap C_H(e)$, $c^q(e)=|C^q(e)|$ and
$c^q_H(e)=|C^q_H(e)|$. Clearly, $c^1(e)+ \cdots + c^s(e)=c(e)$ and $c^1_H(e)+ \cdots + c^s_H(e)=c_H(e)$.
The expectation of $c^q_H(e)$ is $\rho c^q(e)$. Consider some $C^q(e)$ with $c^q(e) > \sqrt{t}$.
According to a large deviation inequality of Chernoff (cf. \cite{AlSp} Appendix A),
for every $\eta > 0$, and in particular for $\eta=\mu/8$, if $n$ (and hence $t$ and hence $c^q(e)$) is sufficiently large,
$$
\Pr[| c^q_H(e) - \rho c^q(e)| > \eta \rho c^q(e) ] < e^{-\frac{2(\eta\rho c^q(e))^2}{c^q(e)}} = e^{-2\eta^2\rho^2c^q(e)} << t^{-k-1}.
$$
It follows that with probability at least $1-st^{-k-1} > 1-t^{-3}$, for all $C^q(e)$ with $c^q(e) > \sqrt{t}$,
$ (1-\eta)\rho c^q(e) \leq c^q_H(e) \leq (1+\eta)\rho c^q(e)$ holds.
Since the sum of $c^q(e)$ having $c^q(e) \le \sqrt{t}$ is  $O(t^{k-2.5})$ and since $c(e)=\Theta(t^{k-2})$ we have
that this sum is much less than $\rho\eta c(e)$. Thus, together with (\ref{e1}) and the fact that $\rho < \psi'(H)^{r-1}\delta^{-r}$
we have
\begin{equation}
\label{e2}
c_H(e) = \sum_{q=1}^s c^q_H(e) \leq \rho(1+\eta)(\sum_{q=1}^s c^q(e)) +\rho \eta c(e) =\rho(1+2\eta)c(e) \leq
\end{equation}
$$
 \rho(1+2\eta)t^{k-2}(\zeta+\frac{\prod_{(s,p) \in E(H_0)}d(s,p)}{d(i,j)})=
(1+2\eta)t^{k-2}(\psi'(H)^{r-1}+\zeta\rho) \leq
$$
$$
t^{k-2}\psi'(H)^{r-1}(1+2\eta)(1+\zeta \delta^{-r})\leq
t^{k-2}\psi'(H)^{r-1}(1+\mu/4)(1+\mu/2) \leq (1+\mu)t^{k-2}\psi'(H)^{r-1}.
$$
Similarly,
\begin{equation}
\label{e3}
c_H(e) \geq \rho(1-\eta)c(e) - \rho \eta c(e) =\rho(1-2\eta)c(e) \geq
\end{equation}
$$
\rho(1-2\eta)t^{k-2}(\frac{\prod_{(s,p) \in E(H_0)}d(s,p)}{d(i,j)}-\zeta)=
(1-2\eta)t^{k-2}(\psi'(H)^{r-1}-\zeta\rho) \geq
$$
$$
t^{k-2}\psi'(H)^{r-1}(1-2\eta)(1-\zeta \delta^{-r})\geq
t^{k-2}\psi'(H)^{r-1}(1-\mu/4)(1-\mu/2) \geq (1-\mu)t^{k-2}\psi'(H)^{r-1}.
$$
Combining (\ref{e2}) and (\ref{e3}) we have that (\ref{e0}) holds for a fixed $e \in E(X_H)$ with probability
at least $1-t^{-3}$. As $E(X_H) < n^2$ we have that (\ref{e0}) holds for all $e \in E(X_H)$ with
probability at least $1-n^2/t^3=1-m^3/n$.
\endpf

We also need the following lemma that gives a lower bound for the number of edges of $X_H$.
\begin{lemma}
\label{l32}
With probability at least $1-1/n$,
$$
|E(X_H)|  > (1-2\zeta)r\frac{n^2}{m^2}\psi'(H).
$$ 
\end{lemma}
{\bf Proof:}\,
We use the notations from Lemma \ref{l31} and the paragraph preceding it.
For $(i,j) \in E(H_0)$, the expected number of edges of $E(V_i,V_j)$ that received the
color $H$ is precisely $d(i,j)\frac{n^2}{m^2}\frac{\psi'(H)}{d(i,j)}=\frac{n^2}{m^2}\psi'(H)$.
Summing over all $r$ edges of $H_0$, the expected number of edges of $W_H$ that received
the color $H$ is precisely $r\frac{n^2}{m^2}\psi'(H)$. As at most $\zeta |E(W_H)|$ edges belong to $W_H$ and
do not belong to $W'_H$ we have that the expectation of $|E(X_H)|$ is at least
$(1-\zeta)r\frac{n^2}{m^2}\psi'(H)$. As $\zeta$, $r$, $m$ are constants and as $\psi'(H)$ is bounded from below by the
constant $m^{1-k_0}$, we have, by the common large deviation inequality of Chernoff (cf. \cite{AlSp} Appendix A), that for
$n > N$ sufficiently large, the probability that $|E(X_H)|$ deviates from its mean by more than 
$\zeta r\frac{n^2}{m^2}\psi'(H)$ is exponentially small in $n$. In particular, the lemma follows. \endpf

Since $R$ contains at most $O(m^{k_0})$ labeled copies of elements of ${\cal F}$ with at most $k_0$ vertices,
we have that with probability at least $1-m^{k_0}/n - m^{k_0+3}/n > 0$ (here we assume again that $N$ is sufficiently large)
{\em all} labeled copies $H$ of elements of ${\cal F}$ in $R$ with $\psi'(H) > m^{1-k_0}$
satisfy the statements of Lemma \ref{l31} and Lemma \ref{l32}.
We therefore fix a coloring for which Lemma \ref{l31} and Lemma \ref{l32} hold
for all labeled copies $H$ of elements of ${\cal F}$ in $R$ having $\psi'(H) > m^{1-k_0}$.

Let $H$ be a labeled copy of some $H_0 \in {\cal F}$ in $R$ with $\psi'(H) > m^{1-k_0}$,
and let $r$ denote the number of edges of $H$.
We construct an $r$-uniform hypergraph $L_H$ as follows.
The vertices of $L_H$ are the edges of the corresponding $X_H$ from Lemma \ref{l31}.
The edges of $L_H$ correspond to the edge sets of the subgraphs of $X_H$ that are partite isomorphic
to $H_0$. We claim that our hypergraph satisfies the conditions of Lemma \ref{l23}.
Indeed, let $q$ denote he number of vertices of $L_H$. Notice that Lemma \ref{l32} provides a lower bound for $q$.
Let $d=t^{k-2} \psi'(H)^{r-1}$. Notice that by Lemma \ref{l31} {\em all} vertices of $L_H$ have their degrees between
$(1-\mu)d$ and $(1+\mu)d$. Also notice that the co-degree of any two vertices of $L_H$ is at most $t^{k-3}$
as two edges cannot belong, together, to more than $t^{k-3}$ subgraphs of $X_H$ that are partite isomorphic to
$H_0$. In particular, for $N$ sufficiently large,
$\mu d > t^{k-3}$. By Lemma \ref{l23} we have at least $(q/r)(1-\beta)$ edge-disjoint copies of $H_0$
in $X_H$. In particular, we have at least
$$
(1-\beta)(1-2\zeta)\frac{n^2}{m^2}\psi'(H) > (1-2\beta)\psi'(H)\frac{n^2}{m^2}
$$
such copies. Recall that $w(\psi')\geq m^2(\alpha-\delta)$. Since there are at most $O(m^{k_0})$
labeled copies $H$ of elements of ${\cal F}$ in $R$ with $0 < \psi'(H) \leq m^{1-{k_0}}$,
their total contribution to $w(\psi')$ is at most $O(m)$.
Hence, summing the last inequality over all $H$ with $\psi'(H) > m^{1-{k_0}}$ we have
at least
$$
(1-2\beta)m^2(\alpha-\delta-O(\frac{1}{m})) \frac{n^2}{m^2} > n^2(\alpha -\epsilon)
$$
edge disjoint copies of elements of ${\cal F}$ in $G$. It follows that $\nu_{{\cal F}}(G) \geq n^2(\alpha-\epsilon)$.
As $\nu^*_{{\cal F}}(G) = \alpha n^2$, Theorem \ref{t1} follows. \endpf

The proof of Theorem \ref{t1} implies an $O(n^{poly(k_0)})$ time algorithm that produces a set of
$n^2(\alpha-\epsilon)$ edge-disjoint copies of elements of ${\cal F}$ in $G$ with probability at least, say, $0.99$.
Indeed, Lemma \ref{l21} can be implemented in $o(n^3)$ time using the algorithm of Alon et. al. \cite{AlDuLeRoYu}.
Lemma \ref{l23} has a polynomial running time implementation due to Grable \cite{Gr}.
Since we only need to compute $\psi^{**}$, rather than $\psi$, we can do this in $O(n^{poly(k_0)})$ time using any polynomial time
algorithm for LP. The other ingredients of the proof are easily implemented in polynomial time.

\end{document}